\documentclass[a4paper]{amsart}
\usepackage{graphics}
\usepackage{color}
\usepackage{amssymb}
\usepackage[all]{xy}
\usepackage{amsmath}
\usepackage{stmaryrd}
\CompileMatrices
\newtheorem*{introthm}{Theorem}

\newtheorem{theorem}{Theorem}[section]
\newtheorem{lemma}[theorem]{Lemma}
\newtheorem{proposition}[theorem]{Proposition}
\newtheorem{corollary}[theorem]{Corollary}
\theoremstyle{definition}
\newtheorem{definition}[theorem]{Definition}
\newtheorem{example}[theorem]{Example}

\newtheorem{remark}[theorem]{Remark}
\theoremstyle{remark}

\def\quot{/\!\!/}
\def\mal{\! \cdot \!}

\def\t#1{\widetilde{#1}}
\def\b#1{\overline{#1}}

\def\CC{{\mathbb C}}
\def\KK{{\mathbb K}}

\def\ZZ{{\mathbb Z}}

\def\QQ{{\mathbb Q}}
\def\PP{{\mathbb P}}

\def\WDiv{\operatorname{WDiv}}

\def\Hom{{\rm Hom}}

\def\Spec{{\rm Spec}}
\def\Proj{{\rm Proj}}
\def\PPDiv{{\rm PPDiv}}
\def\Pol{{\rm Pol}}

\def\cone{{\rm cone}}

\def\eval{{\rm eval}}

\def\vertices{{\rm vertices}}

\newcounter{itemnumber}

\begin{document}
\title[Complete orbit spaces]%
{Complete orbit spaces \\ of affine torus actions}
\author[J.~Hausen]{J\"urgen Hausen} 
\address{Mathematisches Institut, Universit\"at T\"ubingen,
Auf der Morgenstelle 10, 72076 T\"ubingen, Germany}
\email{hausen@mail.mathematik.uni-tuebingen.de}
\begin{abstract}
Given an action of an algebraic 
torus on a normal affine variety,
we describe all open subsets admitting
a complete orbit space.
\end{abstract}

\maketitle

\section*{Introduction}

Let an algebraic torus $T$
act on a normal, algebraic 
variety $X$.
It is an open problem in Geometric 
Invariant Theory to describe 
the collection of all $T$-invariant
open subsets $U \subseteq X$ admitting
a geometric quotient $U \to U/T$
with a complete orbit variety $U/T$.
Several constructions are known
to produce such $U \subseteq X$,
e.g., Mumford's method~\cite{GIT} 
yields in many cases 
subsets $U \subseteq X$ admitting
projective orbit spaces,
and there are more general approaches
providing also non-projective 
complete orbit varieties,
see~\cite{BBSwAmJ}.
However, only in very special cases, 
e.g., $X$ projective and $\dim(T) \le 2$, 
or $X = \PP^n$ or a toric variety,
there are descriptions of {\em all\/} 
$T$-invariant open
subsets $U \subseteq X$
with a complete orbit variety,
see~\cite{BB2},
\cite{BBSo2}, \cite{BBSo3}
and~\cite{BBSw1}.

In the present paper, 
we solve the
above problem for the case 
that an arbitrary torus $T$
acts on a normal, 
{\em affine\/} variety $X$.
Our motivation to consider this 
case is twofold.
Firstly, we hope it to be of use 
for the 
projective and, more generally, the 
divisorial case,
because one can reduce these cases
to the affine one via 
equivariant (multi-)cone 
constructions,
compare~\cite{Ha1} and the 
Example~\ref{grassmann} given
at the end.
Our second motivation concerns
the (in general) non-separated
orbit space $W/T$ of 
the union $W \subseteq X$ 
of all $T$-orbits of maximal dimension.
From a more algebraic point of view, 
$W/T$ is a multigraded analogue of a 
homogeneous spectrum,
compare~\cite{BrS}.
Its complete open subvarieties 
are precisely the complete 
orbit spaces $U/T$,
and thus a description of
them may be helpful for 
a better understanding of $W/T$.

So far, the known approaches
to the affine case basically deal
with diagonal torus actions 
on the affine space $X = \KK^n$.
There are treatments in terms 
of toric geometry, see e.g.~\cite{Hm},
and, alternatively, there is 
a Gale dual approach as presented 
in~\cite{BBSw2}. 
In this paper, we provide 
a general approach, using 
the language of 
proper polyhedral
divisors introduced in~\cite{AlHa}.

A proper polyhedral divisor 
(for short pp-divisor) on a 
normal projective variety $Y$ 
may be written as a linear 
combination of pairwise different
prime divisors
$D_i$ having certain 
polyhedra $\Delta_i$ 
as their coefficients,
which live 
in a common 
rational vector space
and have 
a common pointed tail cone:
\begin{eqnarray*}
\mathfrak{D}
& = & 
\sum_{i=1}^r \Delta_i \otimes D_i.
\end{eqnarray*}
To any such pp-divisor $\mathfrak{D}$
one may associate in a canonical 
way 
a normal affine variety $X$ with
an effective action of a torus $T$,
and, conversely,  
any effective action of a torus 
on a normal affine variety 
is obtained in this way,
see~\cite{AlHa}.
For convenience, we give the precise
definitions and recall the basic
constructions in Section~\ref{sec:ppdiv}.

Given a proper polyhedral
divisor $\mathfrak{D}$ on a projective
variety $Y$ as before,
the basic concept of this paper 
is the notion
of a {\em $\mathfrak{D}$-coherent collection}:
this is a collection of vertices 
$v_i \in \Delta_i$, where $i=1, \ldots, r$,
satisfying certain compatibility
conditions, see Definition~\ref{def:coherent},
which in the case of a curve $Y$ even
turn out to be empty.
The main result is the following,
see Theorem~\ref{mainthm}:

\begin{introthm}
Let $\mathfrak{D}$ be a proper polyhedral
divisor on a normal projective variety~$Y$, 
and let $X$ be the associated normal affine
$T$-variety. 
Then the $\mathfrak{D}$-coherent collections
are in bijection with the 
$T$-invariant open subsets $U \subseteq X$
admitting a geometric quotient $U \to U/T$ 
with a complete orbit space $U/T$.
\end{introthm}

The paper is organized as follows.
In the first section, we recall
among other things 
the language of proper polyhedral 
divisors from~\cite{AlHa},
and we present the basic facts 
needed here.
Section~\ref{sec:prepobs} is devoted 
to preparing investigations
concerning complete orbit spaces.
In Section~\ref{sec:main}, we formulate 
and prove the main result.
Finally, in the last section, we 
discuss an application and examples.

\section{The language of polyhedral divisors}
\label{sec:ppdiv}

In this section, we fix 
(most of) our notation,
give some background on quotients 
and torus actions,
and then recall the necessary 
concepts and results from~\cite{AlHa}.
In particular, we give the precise 
definition of a proper polyhedral 
divisor $\mathfrak{D}$ 
on a semiprojective variety $Y$,
we indicate how to obtain the 
associated affine $T$-variety $X$,
and we describe the fibres of the map 
$\pi \colon \t{X} \to Y$ associated 
to $\mathfrak{D}$.

We work over an algebraically closed 
field $\KK$ of characteristic zero.
By a variety we mean a separated reduced 
$\KK$-scheme of finite type, and the
word prevariety refers to the (possibly)
nonseparated analogue.
By a point of a (pre-)variety, we always 
mean a closed point.

An action $G \times X \to X$
of an algebraic group 
$G$ on a variety $X$ is always  
assumed to be morphical;
in this setting, we also speak 
of the $G$-variety~$X$.
Now suppose that $G$ is reductive, 
for example $G$ is a torus,
and let $X$ be a $G$-variety.
We will have to 
distinguish between the 
following concepts of 
quotients:

\begin{itemize}
\item
A {\em good prequotient\/} 
for the $G$-variety $X$ is 
an affine $G$-invariant 
morphism $\pi \colon X \to Y$ 
onto a (possibly nonseparated) 
prevariety $Y$ such that
$\pi^* \colon \mathcal{O}_Y \to 
\pi_*(\mathcal{O}_X)^G$ is an 
isomorphism.
\item
A {\em geometric prequotient\/} 
for the $G$-variety $X$ is 
a good prequotient $\pi \colon X \to Y$ 
such that each set-theoretical fibre
$\pi^{-1}(y)$, where $y \in Y$, consists
of precisely one $T$-orbit.
\item
A {\em good quotient\/} 
for the $G$-variety $X$ is 
a good prequotient $\pi \colon X \to Y$ 
with a variety $Y$.
\item 
A {\em geometric quotient\/} 
for the $G$-variety $X$ is 
a geometric prequotient 
$\pi \colon X \to Y$ 
with a   variety 
$Y$.
\end{itemize}

If one of these quotients $\pi \colon X \to Y$
exists, then 
it has the following universal property:
let $\varphi \colon X \to Z$ be a 
$G$-invariant morphism to a prevariety
$Z$, then there is a unique morphism
$\psi \colon Y \to Z$ with
$\varphi = \psi \circ \pi$.
This justifies the notations
$Y = X \quot G$ for the good
(pre-)quotient space, 
and $Y = X / G$ in the geometric
case.
We will also refer to $X/G$ as the 
{\em orbit space}.

We shall frequently use two existence
statements on quotients.
Firstly for any affine $G$-variety
$X$, there is a good quotient 
$X \to X \quot G$ with $X \quot G$
being the spectrum of the invariants
$\Gamma(X,\mathcal{O})^G$.
Secondly, if $G$ is a torus, and 
$X$ is a $G$-variety containing only
orbits of maximal dimension,
then there is a geometric prequotient
$X \to X/G$, see~\cite[Corollary~3]{Su}.

\goodbreak

Let us now recall the basic concepts
for actions of algebraic tori $T$
on affine varieties $X$. 
There is a natural correspondence 
between multigraded affine algebras
and such actions:
given a lattice
$M$ and an $M$-graded affine 
algebra
\begin{eqnarray*}
A 
& = & 
\bigoplus_{u \in M} A_u,
\end{eqnarray*}
the torus $T := \Spec(\KK[M])$
acts on the variety 
$X := \Spec(A)$ such that 
the homogeneous elements $f \in A_u$ 
are precisely the semi-invariants
of $X$ with respect to the character
$\chi^u \colon T \to \KK^*$,
and any affine $T$-variety $X$ 
arises in this way.

To the affine $T$-variety $X$ arising from
an $M$-graded affine algebra $A$, we 
may associate combinatorial data 
in terms of $M$. The 
{\em weight cone\/}  
of $X$ is the (convex, polyhedral) 
cone $\omega(X)$
in the rational vector space 
$M_\QQ := \QQ \otimes_\ZZ M$ 
generated by all $u \in M$ 
with $A_u \ne 0$.
The {\em orbit cone\/} of a point
$x \in X$ is the 
(convex, polyhedral) cone $\omega(x)$ 
generated by all $u \in M$ admitting 
an $f \in A_u$ with $f(x) \ne 0$.
Note that the dimension of an 
orbit $\omega(x)$ cone equals 
the dimension of the orbit 
$T \mal x$, and the generic
orbit cone equals the weight 
cone, see~\cite[Section~5]{AlHa}
for a little more background.

We are ready to recall
the necessary
notions and results from~\cite{AlHa}.
In the sequel, $N$ denotes a lattice,
and $\sigma$ is a 
pointed, convex, polyhedral 
cone in the associated rational 
vector space
$N_\QQ = \QQ \otimes_\ZZ N$.
A $\sigma$-polyhedron is a convex 
polyhedron
$\Delta \subseteq N_\QQ$ having
$\sigma$ as its tail cone 
(also referred to as the recession cone).
 
With respect to Minkowski addition, 
the set  $\Pol_{\sigma}^+(N)$ of all 
$\sigma$-polyhedra is a semigroup
with cancellation law; we write 
$\Pol_{\sigma}(N)$ for the associated 
Grothendieck group.
Then the group of {\em polyhedral divisors\/} 
on a normal variety $Y$ is 
$$
\WDiv_{\QQ}(Y,\sigma)
\; := \;
 \Pol_{\sigma}(N) \otimes_{\ZZ} \WDiv_{\QQ}(Y),
$$ 
where $\WDiv_{\QQ}(Y)$ denotes the group
of rational Weil divisors on $Y$.
Via evaluation, any given polyhedral divisor
$\mathfrak{D} = \sum D_i \otimes \Delta_i$ may as 
well be viewed as a piecewise linear convex 
map on the dual cone $\omega \subseteq M_{\QQ}$
of $\sigma \subseteq N_{\QQ}$, 
where $M := \Hom(N,\ZZ)$ is the dual lattice,
namely
$$
\mathfrak{D} \colon \omega \to \WDiv_{\QQ}(Y),
\qquad
u \mapsto \sum \eval_u(\Delta_i) D_i,
\quad
\text{where }
\eval_u(\Delta_i) \; := \; \min_{v \in \Delta_i} u(v).
$$
Here, convexity has to be understood in 
the setting of divisors, 
that means that we always
have  
$\mathfrak{D}(u+u') \ge \mathfrak{D}(u)+\mathfrak{D}(u')$.
A {\em proper polyhedral divisor\/}
(abbreviated pp-divisor)
is a polyhedral divisor 
$\mathfrak{D} \in \WDiv_{\QQ}(Y,\sigma)$
such that
\begin{enumerate}
\item 
there is a representation
$\mathfrak{D} = \sum D_i \otimes \Delta_i$ with
effective $D_i \in \WDiv_\QQ(Y)$ and 
$\Delta_i \in \Pol^{+}_{\sigma}(N)$,
\item 
each evaluation $\mathfrak{D}(u)$, where $u \in \omega$,
is a semiample $\QQ$-Cartier divisor, 
i.e., has a base point free multiple,
\item 
for any $u$ in the relative interior 
$\omega^\circ \subseteq \omega$,
the some multiple of $\mathfrak{D}(u)$ 
is a big divisor, i.e., admits a section
with affine complement.
\end{enumerate}

Now suppose that $Y$ is semiprojective,
i.e., projective over some affine variety,
and let $\mathfrak{D} = \sum D_i \otimes \Delta_i$ 
be a pp-divisor on~$Y$.
Then $\mathfrak{D}$ defines a sheaf
of $\mathcal{O}_Y$-algebras, and we have the 
corresponding relative spectrum: 
$$
\mathcal{A} 
\; := \;
\bigoplus_{u \in \omega \cap M}
\mathcal{O}(\mathfrak{D}(u)),
\qquad
\t{X} 
\; := \;
\Spec_Y(\mathcal{A}).
$$
The grading of $\mathcal{A}$ gives rise to 
an effective action of the torus $T := \Spec(\KK[M])$
on $\t{X}$, and the canonical map 
$\pi \colon \t{X} \to Y$ is a good quotient 
for this action.

By~\cite[Theorem~3.1]{AlHa},
the ring of global sections
$A := \Gamma(\t{X},\mathcal{O}) = \Gamma(Y, \mathcal{A})$
is finitely generated and normal, 
and there is a $T$-equivariant, birational,
proper morphism $r \colon \t{X} \to X$
onto the normal, affine $T$-variety
$X := \Spec(A)$.
Conversely~\cite[Theorem~3.4]{AlHa},
says that every normal, affine variety
with an effective torus action arises 
in the above way from a pp-divisor on a
semiprojective variety.

\begin{remark}
\label{goodaction}
For the affine $T$-variety $X$ 
arising from a pp-divisor
$\mathfrak{D}$ 
on a semiprojective variety 
$Y$, the following 
statements are equivalent:
\begin{enumerate}
\item
All $T$-orbits of $X$ have 
a common orbit $T \mal x_0$
in their closures.
\item
The weight cone $\omega(X)$ is pointed, 
and $A_0 = \KK$ holds. 
\item
The semiprojective variety $Y$ is projective.
\end{enumerate}
\end{remark}

We will need parts of the 
description of the fibres of the 
map
$\pi \colon \t{X} \to Y$ given 
in~\cite[Prop.~7.8 and Cor.~7.9]{AlHa}.
First recall that,
for a $\sigma$-polyhedron 
$\Delta$ in $N_{\QQ}$,
each face $F \preceq \Delta$
defines a convex, polyhedral 
cone in $M_\QQ$ via
\begin{eqnarray*}
F 
& \mapsto &
\lambda(F)
\; := \; 
\{u \in M_{\QQ}; \; 
\langle u, v - v' \rangle \ge 0
\text{ for all } 
v \in \Delta, \, v' \in F\}.
\end{eqnarray*}
The collection $\Lambda(\Delta)$ 
of all these cones is called the normal 
quasifan of $\Delta$;
it subdivides the dual cone
$\omega \subseteq M_\QQ$
of $\sigma \subseteq N_\QQ$.
Note that the
normal quasifan 
$\Lambda(\Delta_1 + \Delta_2)$
of a Minkowski sum is 
the coarsest common refinement
of $\Lambda(\Delta_1)$
and $\Lambda(\Delta_2)$.

Now, let 
$\mathfrak{D} = \sum \Delta_i \otimes D_i$
be a representation of our pp-divisor
such that all $D_i$ are prime.
For a point $y \in Y$, its 
{\em fiber polyhedron\/} is 
the Minkowski sum
$$
\Delta_y 
\; := \;
\sum_{y \in D_i} \Delta_i
\; \in \; 
\Pol_{\sigma}^+(N).
$$

\begin{theorem}
\label{thm:fibres}
Let $y \in Y$, consider the
affine $T$-variety $\pi^{-1}(y)$,
and let $\Lambda_y$ denote the 
normal quasifan of the fiber
polyhedron $\Delta_y$. 
Then there is a one-to-one correspondence: 
$$
\{ 
T \text{-orbits in } \pi^{-1}(y)
\}
\;  \to \;
\Lambda_y
\qquad
T \mal \t{x}
\; \mapsto \;
\omega(\t{x}).
$$
\end{theorem}

Secondly, we shall need
parts of the description 
of the $T$-orbits
of $X$ given in~\cite[Theorem~10.1]{AlHa}.
This
involves the canonical contraction
maps
$$ 
\vartheta_u
\colon 
Y 
\; \to \; 
\Proj\left( 
\bigoplus_{n=0}^{\infty}
\Gamma(Y,\mathcal{O}(\mathfrak{D}(nu)))
\right),
\quad
\text{where } u \in \omega \cap M.
$$

\begin{theorem}
\label{thm:orbits}
For any two  
$\t{x}_1, \t{x}_2 \in \t{X}$,
the following statements
are equivalent:
\begin{enumerate}
\item
The contraction morphism $r \colon \t{X} \to X$
identifies the orbits $T \mal \t{x}_1$ and $T \mal \t{x}_2$.
\item
We have $\omega(\t{x}_1) = \omega(\t{x}_2)$ and 
$\vartheta_u(\pi(\t{x}_1)) = \vartheta_u(\pi(\t{x}_2))$
for some $u \in \omega(\t{x}_1)^\circ$.
\end{enumerate}
\end{theorem}

\section{Preparing observations}
\label{sec:prepobs}

In this section,
$X$ is the normal, 
affine $T$-variety arising 
from a pp-divisor $\mathfrak{D}$
living on a normal,
semiprojective variety~$Y$.
As before, 
$r \colon \t{X} \to X$
denotes 
the associated $T$-equivariant 
birational contraction map,
and $\pi \colon \t{X} \to Y$ 
is the associated 
good quotient for the $T$-action.

We show that existence of a complete 
orbit space $U/T$ for a subset 
$U \subseteq X$ is equivalent
to existence of a complete orbit 
space $\t{U}/T$ for 
$\t{U} := r^{-1}(U)$,
and we give a geometric
characterization of the 
subsets $\t{U} \subseteq \t{X}$
admitting a complete orbit space.
We establish these facts in a 
series of Lemmas, and then 
gather them in 
Proposition~\ref{prop:simple}.

\begin{lemma}
\label{categorical}
Let $\t{U} \subseteq \t{X}$ 
be a $T$-invariant
open subset.
Then $Y' := \pi(\t{U})$
is open in~$Y$, and,
for any $T$-invariant morphism 
$\varphi \colon \t{U} \to Z$
to a   variety $Z$,
there is a unique morphism 
$\psi \colon Y' \to Z$ with
$\varphi = \psi \circ \pi$.
\end{lemma}

\begin{proof}
We first consider any affine 
open subset $Y_0 \subseteq Y$. 
Then also $\t{X}_0 := \pi^{-1}(Y_0)$
is affine, and hence 
$\t{U}_0 := \t{U} \cap \t{X}_0$ is 
a union of homogeneous 
localizations $\t{U}_f := (\t{X}_0)_f$.
For each of these localizations,
we have a commutative diagram
$$ 
\xymatrix{
{\t{U}_f} 
\ar[r]
\ar[d]_{\pi_f}^{\quot T}
&
{\t{X}_0}
\ar[d]^{\pi}_{\quot T}
\\
{\t{U}_f \quot T} 
\ar[r]_{\imath_f}
&
Y_0
}
$$

Using e.g.~Theorem~\ref{thm:fibres},
we see that the generic fiber of 
$\pi \colon \t{X}_0 \to Y_0$
is the closure of a single 
$T$-orbit.
The above diagram
tells us that the same must hold
for the quotient map
$\pi_f \colon \t{U}_f \to \t{U}_f \quot T$.
Consequently, 
we have canonical isomorphisms
$$ 
\KK( \t{U}_f \quot T) 
\; \cong \; 
\KK( \t{U}_f)^T
\; = \; 
\KK(\t{X}_0)^T 
\; \cong \;
\KK(Y_0).
$$
Thus, 
$\imath_f \colon \t{U}_f \quot T \to Y_0$
is birational.
By 
Theorem~\ref{thm:fibres},
the fibers of
$\pi$ contain only finitely many 
$T$-orbits.
Thus, $\imath_f$ 
has finite fibers, and 
hence is an open embedding.
Since $Y' \cap Y_0$ is covered by 
the images 
$Y_f := \imath_f(\t{U}_f \quot T)$,
it must be open in $Y_0$.

Given a $T$-invariant morphism
$\varphi_0 \colon \t{U}_0 \to Z$
to a variety $Z$,
consider any restriction
$\varphi_f \colon \t{U}_f \to Z$.
The above consideration yields 
a unique morphism
$\psi_f \colon Y_f \to Z$ 
with 
$\varphi_f = \psi_f \circ \pi$.
Moreover, any two 
such $\psi_f, \psi_g$
coincide on the dense
subset 
$$
Y_{fg}
\; = \; 
\pi(\t{U}_f \cap \t{U}_g)
\; \subseteq \; 
Y_f \cap Y_g.
$$
Consequently,
since $Z$ is separated, 
we can glue together the 
morphisms $\psi_f \colon Y_f \to Z$ 
to a morphism 
$\psi_0 \colon Y' \cap Y_0 \to Z$,
and obtain this way a unique 
factorization
$\varphi_0 = \psi_0 \circ \pi$.

To conclude the proof, 
cover $Y$ by affine open 
subsets $Y_i$.
Then, by the preceding 
consideration, 
each $Y_i' := Y' \cap Y_i$ is open,
and hence 
$Y' \subseteq Y$ is so. 
Moreover, given a $T$-invariant 
$\psi \colon \t{U} \to Z$ to a 
variety $Z$,
we have a factorization 
$\varphi = \psi_i \circ \pi$ 
over each $Y_i'$, and, 
by uniqueness over $Y_i' \cap Y_j'$,
the $\psi_i$ can be patched together to 
the desired morphism
$\psi \colon Y' \to Z$.
\end{proof}

\begin{lemma}
\label{sepfibers}
Let $\t{U} \subseteq \t{X}$
be a $T$-invariant, open
subset
containing only $T$-orbits
of maximal dimension,
and set $Y' := \pi(U)$.
Then the following
statements are equivalent:
\begin{enumerate}
\item 
The orbit space $\t{U} / T$
is separated.
\item 
$\t{U} \cap \pi^{-1}(y)$
is a single $T$-orbit for every $y \in Y'$.
\end{enumerate}
If one of these statements holds, 
then  the restriction 
$\pi \colon \t{U} \to Y'$
is a geometric quotient for 
the $T$-action.
\end{lemma}

\begin{proof}
Recall from Lemma~\ref{categorical}
that $Y' = \pi(\t{U})$ is open in
$Y$.
Suppose that~(i) holds.
Then 
Lemma~\ref{categorical}
and the universal property of
$\t{U} \to \t{U} / T$ yield
that the canonical morphism 
$\t{U} /T \to Y'$ is an isomorphism.
In particular, $\pi \colon \t{U} \to Y'$
is a geometric quotient, 
and~(ii) holds.

Suppose that~(ii) holds.
Then it suffices to show that 
$\pi \colon \t{U} \to Y'$
is a geometric quotient.
First we note that 
$\t{U}$ can be covered
by $T$-invariant open affine
subsets $\t{U}_0 \subseteq \t{U}$,
see~\cite[Cor.~2]{Su}.
For each such $\t{U}_0$,
we obtain a commutative diagram
$$ 
\xymatrix{
{\t{U}_0}
\ar[r]^{\subseteq}
\ar[d]_{/T}
&
{\t{U}}
\ar[d]^{\pi}
\\
{\t{U}_0 / T} 
\ar[r]_{\imath}
&
{Y'}
}
$$
The induced map 
$\imath \colon \t{U}_0 / T \to Y'$
is birational, and, by assumption,
injective.
Hence it is an open embedding,
and $\pi(\t{U}_0)$ is affine.
Thus, 
$\pi \colon \t{U} \to Y'$ 
looks locally w.r. to $Y'$ 
like an affine 
geometric quotient,
and hence is a geometric 
quotient.
\end{proof}

\begin{lemma}
\label{lem:compl}
Let $U \subseteq X$ be a $T$-invariant 
open subset containing only $T$-orbits 
of maximal dimension, and set
$\t{U} := r^{-1}(U)$.
Then the following statements are 
equivalent:
\begin{enumerate}
\item
The orbit space $U/T$ is a complete
  variety.
\item
The orbit space $\t{U}/T$ is a 
complete   variety.
\end{enumerate}
In each of these two cases, 
$Y = \pi(\t{U})$ holds, 
$Y$ is projective, 
and $\pi \colon \t{U} \to Y$
is a geometric quotient;
in particular, $\t{U}/T$ 
is then projective.
\end{lemma}

\begin{proof}
Note that $\t{U} = r^{-1}(U)$ contains only 
orbits of maximal dimension.
Thus, there is a geometric 
prequotient $\t{U} \to \t{U} / T$,
and we have a commutative diagram
$$ 
\xymatrix{
{\t{U}}
\ar[r]^{r}
\ar[d]_{/T}
&
{U}
\ar[d]^{/T}
\\
{\t{U} / T} 
\ar[r]_{\imath}
&
{U / T}
}
$$

If~(ii) holds, then 
we may apply~\cite[Lemma~3.2]{BBSw1}
to the (birational) surjective morphism 
$\t{U}/T \to U /T$,
and obtain that $U/T$ is a complete 
variety. 

Now suppose that~(i) holds.
If $\t{U}/T$ is not separated, 
then Lemma~\ref{sepfibers} provides 
two different orbits  
$T \mal \t{x}_1$ and $T \mal \t{x}_2$
in $\t{U}$, which lie 
in a common fibre 
$\pi^{-1}(y) \subseteq \t{X}$.
By Lemma~\ref{categorical},
their images $y_i \in \t{U} / T$
are identified to a point $y \in U /T$
under
$\imath \colon  \t{U} / T \to U / T$. 
Let $x \in U$ lie over  $y \in U /T$.
Then $r \colon \t{X} \to X$ maps 
each orbit $T \mal \t{x}_i$ onto 
$T \mal x$.
By Theorems~\ref{thm:fibres} 
and~\ref{thm:orbits},
this is impossible for two different
$T$-orbits inside one fibre 
$\pi^{-1}(y) \subseteq \t{X}$.
Thus,  
$\t{U} / T$ must be separated.

In order to see that $\t{U} / T$
is complete, it suffices to show
that $\t{U} / T \to U/T$ is a 
proper morphism.
Since $\t{U}/T$ is a variety,
$\t{U}/T \to U/T$ is of finite type
and separated.
Universal closedness follows 
directly from that fact that,
given any morphism 
$Z \to U/T$, we have a 
canonical commutative
diagram
$$ 
\xymatrix{
Z \times_{U/T} \t{U}
\ar[rr]^{\rm proper}
\ar[d]_{/T}
&
&
Z \times_{U/T} U
\ar[d]^{/T}
\\
Z \times_{U/T} \t{U} /T
\ar[rr] \ar[dr] \ar[ddr]
&
&
Z \times_{U/T} U /T
\ar[dl]_{\cong}
\ar[ddl]
\\
& Z \ar[d] &
\\
&
U/T
&
}
$$

Knowing that $\t{U}/T$ is a complete variety,
we can conclude that the canonical
(dominant) morphism $\t{U}/T \to Y$ is 
surjective, which implies $\pi(\t{U}) = Y$.
Lemma~\ref{categorical} then even says that 
$\t{U}/T \to Y$ is an isomorphism.
In particular, 
$\t{U}/T$ is projective and
$\pi \colon \t{U} \to Y$
is a geometric quotient.
\end{proof}

\begin{corollary}
\label{complchar}
There exists an open, $T$-invariant 
subset $U \subseteq X$ admitting a
complete orbit variety $U/T$ if and 
only if $Y$ is projective.
\end{corollary}

As announced before, we now gather the 
observations made in the preceding 
Lemmas. For this, we introduce the 
following notion. 

\begin{definition}
We say that an open subset 
$\t{U} \subseteq \t{X}$ 
is {\em simple\/} if
$\pi(\t{U}) = Y$ holds,
we have $r^{-1}(r(\t{U})) = \t{U}$,
and for every $y \in Y$
the set $\pi^{-1}(y) \cap \t{U}$ 
is a single $T$-orbit.
\end{definition}

\begin{proposition}
\label{prop:simple}
Let $X$ the affine 
$T$-variety arising from
a pp-divisor living on a 
projective variety $Y$.
Then the assignments 
$\t{U} \mapsto r(\t{U})$ 
and $U \mapsto r^{-1}(U)$
define mutually inverse
one-to-one correspondences
between the simple 
subsets of $\t{U} \subseteq \t{X}$ 
and the $T$-invariant
open subsets $U \subseteq X$ 
with a complete
orbit space $U/T$.
\end{proposition}

\section{Complete orbit spaces}
\label{sec:main}

In this section, we formulate and prove 
our main result describing the
open subsets with a complete orbit 
space for a given normal affine 
variety $X$ with an effective 
torus action $T \times X \to X$.
According to 
Corollary~\ref{complchar},
we may assume that 
the $T$-variety
$X$ arises from a pp-divisor on a 
projective variety $Y$;
characterizations of this 
case were given in 
Remark~\ref{goodaction}.

Here comes the precise 
setup of this section.
By $Y$ we denote, as indicated,
a normal, projective variety, 
$N$ is a lattice
and $\sigma \subseteq N_\QQ$ 
is a pointed cone.
Let $\mathfrak{D} \in \PPDiv(Y,\sigma)$
be a pp-divisor  on $Y$, given
by a representation
\begin{eqnarray}
\label{eq:ppdiv}
\mathfrak{D}
& = & 
\sum_{i=1}^r \Delta_i \otimes D_i
\end{eqnarray}
with pairwise different
prime divisors $D_i \in \WDiv(Y)$ 
and $\sigma$-polyhedra 
$\Delta_i \subseteq N_\QQ$.
As before, we denote by $X$ the 
normal, affine $T$-variety arising
from $\mathfrak{D}$, by 
$r \colon \t{X} \to X$ the 
$T$-equivariant contraction map
and by $\pi \colon \t{X} \to Y$ 
the associated good quotient.

Recall from Section~\ref{sec:ppdiv}
that for any $y \in Y$, there is 
an associated fiber polyhedron 
$\Delta_y \subseteq N_\QQ$, and
the normal quasifan $\Lambda_y$ 
of $\Delta_y$
subdivides the dual cone 
$\omega \subseteq M_\QQ$
of $\sigma \subseteq N_\QQ$.
We have the bijection $F \mapsto \lambda(F)$
from the faces $F \preceq \Delta_y$ 
to the cones of~$\Lambda_y$.
For a cone $\lambda \subseteq M_\QQ$,
we denote its relative interior
 by $\lambda^{\circ}$.

Let us introduce the combinatorial 
data for the description of the 
collection of all $T$-invariant,
open subsets $U \subseteq X$ 
admitting a complete orbit
variety $U/T$. The definition
makes use of the canonical contraction
maps mentioned in Section~\ref{sec:ppdiv}:
$$ 
\vartheta_u
\colon 
Y 
\; \to \; 
\Proj\left( 
\bigoplus_{n=0}^{\infty}
\Gamma(Y,\mathcal{O}(\mathfrak{D}(nu)))
\right),
\quad
\text{where } u \in \omega \cap M.
$$

\begin{definition}
\label{def:coherent}
Let 
$\mathfrak{D} = \sum_{i=1}^r \Delta_i \otimes D_i$
be a pp-divisor on a normal, projective variety~$Y$
as in~(\ref{eq:ppdiv}),
and consider vertices $v_i \in \Delta_i$, 
where $i = 1, \ldots, r$.
We say that $v_1, \ldots, v_r$ 
is a {\em $\mathfrak{D}$-admissible collection\/} 
if for any $y \in Y$ the point
\begin{eqnarray*}
v_y 
& := 
\sum_{y \in D_i} v_i
\end{eqnarray*} 
is a vertex of $\Delta_y$.
If $v_1, \ldots, v_r$ is a 
$\mathfrak{D}$-admissible collection
and $y \in Y$ is given, we write
$\lambda_y := \lambda(v_y) \in \Lambda_y$ 
for the corresponding cone. 
We say that a $\mathfrak{D}$-admissible collection 
$v_1, \ldots, v_r \in N_\QQ$ is
{\em $\mathfrak{D}$-coherent\/} if for any two
$y_1,y_2 \in Y$ we have
\begin{eqnarray*}
\lambda_{y_2} \in \Lambda_{y_1}
\text{ and }
\vartheta_u(y_2) = \vartheta_u(y_1)
\text{ for some }
u \in \lambda_{y_2}^\circ
& \implies &
\lambda_{y_1}
= \lambda_{y_2}.
\end{eqnarray*}
\end{definition}

Note that, if all divisors $\mathfrak{D}(u)$ 
corresponding to interior points of $\omega$ 
are ample, then their contraction maps $\vartheta_u$ 
are trivial, and thus, 
every $\mathfrak{D}$-admissible collection
is $\mathfrak{D}$-coherent. 
This holds
for example if $Y$ is a curve, 
or if $Y$ has a free cyclic 
divisor class group.

\begin{definition}
Let 
$\mathfrak{D} = \sum_{i=1}^r \Delta_i \otimes D_i$
be a pp-divisor on a normal, 
projective variety~$Y$ as 
in~(\ref{eq:ppdiv}).
To any $\mathfrak{D}$-admissible collection
$v_1, \ldots, v_r$ we associate $T$-invariant
subsets
\begin{eqnarray*}
\t{U}(v_1, \ldots, v_r)
& := &
\{\t{x} \in \t{X}; \; 
\omega(\t{x}) = \lambda_{\pi(\t{x})} \}
\\
& \subseteq &
\t{X},
\\
U(v_1, \ldots, v_r)
& := &
r(\t{U}(v_1, \ldots, v_r))
\\
& \subseteq &
X.
\end{eqnarray*}
\end{definition}

\begin{theorem}
\label{mainthm}
Let $\mathfrak{D} = \sum_{i=1}^r \Delta_i \otimes D_i$
be a pp-divisor on a normal, projective variety~$Y$
as in~(\ref{eq:ppdiv}), 
and let $X$ be the associated normal, affine 
$T$-variety.
Then there is a bijection:
\begin{eqnarray*}
\{\mathfrak{D} \text{-coherent collections}\}
& \to & 
\left\{
\vcenter{
\hbox{$T$-invariant open $U \subseteq X$ with}
\hbox{a complete orbit space $U/T$}
}
\right\}
\\
(v_1, \ldots, v_r)
& \mapsto &
U(v_1, \ldots, v_r).
\end{eqnarray*}
\end{theorem}

Using the descriptions of projective orbit spaces and 
such embeddable into some toric variety in terms of 
orbit cones given in~\cite[Sec.~1]{ArHa}, we can easily 
figure out such orbit spaces from their defining coherent 
collections, provided that the $T$-variety $X$ is 
factorial.

\begin{remark}
Let $\mathfrak{D}= \sum_{i=1}^r \Delta_i \otimes D_i$ 
be a pp-divisor on a normal, 
projective variety~$Y$, and suppose that the associated  
normal, affine $T$-variety $X$ is factorial.
Let $v_1, \ldots, v_r$ be a coherent collection,
and denote by $U = U(v_1, \ldots, v_r) \subseteq X$ the 
associated open set with complete orbit space. 
\begin{enumerate}
\item
The orbit space $U / T$ is projective if and only if 
the intersection of all relative interiors $\lambda_y^\circ$,
where $y \in Y$ is nonempty.
\item
The orbit space $U / T$ admits an embedding into a toric
variety if and only if for any two $y_1,y_2 \in Y$,
the intersection $\lambda_{y_1}^\circ \cap \lambda_{y_2}^\circ$
is nonempty.
\end{enumerate}
\end{remark}

We turn to the proof of Theorem~\ref{mainthm}.
We shall make use of the 
following elementary observation 
in convex geometry, which is
evident from the definition
of the normal quasi-fan of
a polyhedron.

\begin{lemma}
\label{lem:minkowski}
Let $\Delta_1, \ldots, \Delta_r$ be
polyhedra in a common vector space,
and fix vertices 
$v_i \in \Delta_i$.
Denote by $\Lambda_i := \Lambda(\Delta_i)$ 
the normal quasifan, 
and let $\lambda_i \in \Lambda_i$
be the cone corresponding to $v_i$.
Then the following statements are equivalent.
\begin{enumerate}
\item
The point
$v := v_1 + \ldots + v_r$ is a vertex of 
the Minkowski sum
$\Delta := \Delta_1 + \ldots + \Delta_r$.
\item 
The cone
$\lambda := \lambda_1 \cap \ldots \cap \lambda_r$ 
is a maximal cone of the normal quasifan 
$\Lambda := \Lambda(\Delta)$.
\end{enumerate}
If one of these statements holds, then
$\lambda \in \Lambda$ is the maximal cone 
corresponding to the vertex $v \in \Delta$.
\end{lemma}

\begin{proof}[Proof of Theorem~\ref{mainthm}]
According to Proposition~\ref{prop:simple},
it suffices to show that the assignment 
$(v_1, \ldots, v_r) \mapsto \t{U}(v_1, \ldots, v_r)$
defines a bijection from
the $\mathfrak{D}$-coherent 
collections to 
the simple subsets
$\t{U} \subseteq \t{X}$.

Let $v_1, \ldots, v_r$ be a 
$\mathfrak{D}$-coherent collection. 
Our first task is to check
that the $T$-invariant subset
$\t{U}(v_1, \ldots, v_r) \subseteq \t{X}$
is indeed open.
For this, let $\t{x}_0 \in \t{U}(v_1, \ldots, v_r)$,
and set $y_0 := \pi(\t{x}_0)$.
Then $y_0$ admits a 
canonical open neighbourhood 
$$
V
\; := \;
Y \setminus \bigcup_{y_0 \not \in D_j} D_j
\; \subseteq \;  
Y.
$$
Then, for every $y \in V$, 
the normal quasifan $\Lambda_{y_0}$ of 
$\Delta_{y_0}$ refines the normal
quasifan $\Lambda_y$ of $\Delta_y$, 
and, by Lemma~\ref{lem:minkowski},
for the relative interiors 
$\lambda_{y_0}^\circ \subseteq \lambda_{y_0}$ 
and 
$\lambda_y^\circ \subseteq \lambda_y$ 
of the cones corresponding to 
the vertices 
$v_{y_0} \preceq \Delta_{y_0}$
and 
$v_y \preceq \Delta_y$
we have 
\begin{eqnarray}
\label{eqn:conescont}
\lambda_{y_0}^\circ 
& \subseteq & \lambda_y^\circ.
\end{eqnarray}

Choose an affine open neighbourhood
$V_0 \subseteq V$ of $y_0 \in V$, 
and an integral vector 
$u \in \lambda_{y_0}^\circ$
admitting a homogeneous function
$f \in \Gamma(\pi^{-1}(V_0),\mathcal{O})_u$
with $f(\t{x}_0) \ne 0$.
This gives an open neighbourhood
of $\t{x}_0$ in $\t{X}$ , namely
\begin{eqnarray*}
\pi^{-1}(V_0)_f
& \subseteq &
\{\t{x} \in \pi^{-1}(V_0); \; u \in \omega(\t{x}) \}.
\end{eqnarray*}
According to (\ref{eqn:conescont}), 
the whole set on right hand side 
is contained in $\t{U}(v_1, \ldots, v_r)$.
This implies openness of the subset
$\t{U}(v_1, \ldots, v_r) \subseteq \t{X}$.

Now have to verify the properties of 
a simple set for $\t{U}(v_1, \ldots, v_r)$.
By Theorem~\ref{thm:fibres},
the image 
$\pi(\t{U}(v_1, \ldots, v_r))$
equals $Y$, 
and 
each fibre $\pi^{-1}(y)$, 
where $y \in Y$ contains exactly 
one $T$-orbit
of $\t{U}(v_1, \ldots, v_r)$.
So, we only have to show 
that $\t{U}(v_1, \ldots, v_r)$
is saturated with respect 
to the contraction map
$r \colon \t{X} \to X$.

For this, 
let $\t{x}_1 \in \t{U}(v_1, \ldots, v_r)$
and $\t{x}_2 \in \t{X}$
with $r(\t{x}_2) = r(\t{x}_1)$.
Set $y_i := \pi(\t{x}_i)$.
Then, by Theorem~\ref{thm:orbits}, we have 
$\omega(\t{x}_2) = \omega(\t{x}_1)$
and $\vartheta_u(y_2) = \vartheta_u(y_1)$
for some $u \in \omega(\t{x}_2)^\circ$.
This implies 
$\t{x}_2 \in \t{U}(v_1, \ldots, v_r)$,
because by $\mathfrak{D}$-coherence 
of the collection $v_1, \ldots, v_r$,
we have
$$ 
\omega(\t{x}_2)
\; = \; 
\omega(\t{x}_1)
\; = \; 
\lambda_{y_1}
\; = \; 
\lambda_{y_2}.
$$

Now, let $\t{U} \subseteq \t{X}$ be any
simple subset.
We have to show that $\t{U}$ arises 
from a $\mathfrak{D}$-coherent 
collection.
Recall from Lemma~\ref{sepfibers}
that the restriction
$\pi \colon \t{U} \to Y$ is a 
geometric quotient.
Moreover, we have 
the prime divisors $D_i$ 
in $Y_i$,
and (nonempty) locally 
closed subsets
$$
Y_i  
\; := \;
D_i \setminus \bigcup_{j \ne i} D_j
\; \subseteq \; 
Y,
\qquad
\t{U}_i
\; := \; 
\t{U} \cap \pi^{-1}(Y_i)
\; \subseteq \; 
\t{X}.
$$

Note that 
$\pi \colon \t{U}_i \to Y_i$ 
is a geometric quotient, and hence
$\t{U}_i$ is irreducible.
Moreover, all points $y \in Y_i$
have the same fiber polyhedron 
$\Delta_y = \Delta_i$,
and thus, since $\t{U}_i$ is 
irreducible, 
$\omega(\t{x})$
is constant along $\t{U}_i$.
Finally, since also 
$\t{U} \cap \pi^{-1}(D_i)$ 
is irreducible,
the closure of $\t{U}_i$ 
in $\t{U}$ is given by 
$$
E_i
\; := \; 
\b{\t{U}_i}
\; = \;
\t{U} \cap \pi^{-1}(D_i).
$$

For $\t{x} \in \t{U}_i$,
set 
$\lambda_i := \omega(\t{x})$.
Theorem~\ref{thm:fibres}
tells us that $\lambda_i$ is 
a maximal cone of 
the normal quasifan $\Lambda_i$
of $\Delta_i$.
Let $v_i \in \Delta_i$ denote the 
vertex corresponding to 
$\lambda_i \in \Lambda_i$.
We show that 
$v_1, \ldots, v_r$ is 
a $\mathfrak{D}$-admissible
collection with 
$\t{U} = \t{U}(v_1, \ldots, v_r)$.
By 
Lemma~\ref{lem:minkowski}
it suffices to show that
for any subset 
$I \subseteq \{1, \ldots, r\}$
we have
\begin{equation}
\label{eqn:oc}
\omega(\t{x})
\; = \;
\lambda_I 
\; := \; 
\bigcap_{i \in I} \lambda_i
\quad
\text{for every } 
\t{x} 
\; \in  \; 
\t{U}_I
\; := \; 
\bigcap_{i \in I} E_i
\setminus 
\bigcup_{j \not\in I} E_j.
\end{equation}

Since $\t{U}$ is a union
of subsets $r^{-1}(X_f)$ 
with homogeneous 
$f \in \Gamma(X,\mathcal{O})$,
and $\t{U}_I$ is contained in the
closure of each $\t{U}_i$ with 
$i \in I$, we must have
$\omega(\t{x}) \subseteq \lambda_i$
for all $\t{x} \in \t{U}_I$ and 
all $i \in I$.
Moreover, in this situation,
the fiber polyhedron
$\Delta$ of $\pi(\t{x})$
is the Minkowski
sum of the $\Delta_i$,
where $i \in I$.
The normal quasifan 
$\Lambda := \Lambda(\Delta)$
is the coarsest common refinement 
of the normal quasifans
$\Lambda_i = \Lambda(\Delta_i)$.
Thus, $\lambda_I$ is a maximal
cone of $\Lambda$.
Since $\omega(\t{x}) \in \Lambda$
holds and $\omega(\t{x})$ is of 
full dimension,
we obtain 
$\omega(\t{x}) = \lambda_I$.
 
Finally, we have to verify 
$\mathfrak{D}$-coherence 
of the $\mathfrak{D}$-admissible collection
$v_1, \ldots, v_r$.
So, let $y_1, y_2 \in Y$ such that
$\lambda_{y_2} \in \Lambda_{y_1}$
holds and we have
$\vartheta_u(y_2)=\vartheta_u(y_1)$
for some integral 
$u \in \lambda_{y_2}^\circ$.
By Theorem~\ref{thm:orbits},
the $T$-orbits 
$T \mal \t{x}_i \in \pi^{-1}(y_i)$
corresponding to $\lambda_{y_2}$
are identified under 
$r \colon \t{X} \to X$.
Since $\t{U}$ is saturated w.r. 
to $r \colon \t{X} \to X$, 
we obtain $\t{x}_1 \in \t{U}$,
which implies
$\lambda_{y_1} = \lambda_{y_2}$.
\end{proof}

\section{Applications and Examples}

In this section, we discuss a few 
examples and applications.
The first observation concerns the 
limit $Y'$ over all GIT-quotients
associated to possible 
linearizations of the trivial 
bundle. 
The limit $Y'$ contains a canonical 
component $Y'_0$ compare dominating 
all the GIT-quotients just mentioned,
see e.g~\cite[Section~6]{AlHa}.
We obtain that the normalization
$Y$ of $Y'_0$
dominates moreover all complete
orbit spaces, i.e., also those
that do not arise from GIT: 

\begin{corollary}
Let $U \subseteq X$ admit a
complete orbit space $U(T)$.
Then there is a surjective 
birational 
morphism $Y \to U/T$ 
from the normalized canonical 
component $Y$ of the limit 
over all GIT-quotients of 
$X$.
\end{corollary}

\begin{proof}
The $T$-variety $X$ admits 
a description by a pp-divisor
living $\mathfrak{D}$ living 
on the normalized canonical 
component $Y$,
see~\cite[Section~6]{AlHa}.
Thus, the claim follows from
Lemma~\ref{lem:compl}.
\end{proof}

We now use our result to
treat an example of 
A.~Bia\l ynicki-Birula
and 
J.~\'Swi\c{e}cicka
of a $\KK^*$-action
on the Grassmannian $G(2,4)$, 
see~\cite{BBSw1};
to our knowledge, this the 
simplest example admitting
complete orbit spaces that
are not embeddable into any
toric variety.

\begin{example}
\label{grassmann}
We consider the cone $X$ over the
Grassmannian $G(2,4)$.
In terms of Pl\"ucker
Coordinates, $X$ is given as
$$ 
X 
\; = \; 
V(\KK^6,z_1z_6 - z_2z_5 + z_3z_4)
\; \subset \; \KK^6
$$
Let the twodimensional torus
$T :=  \KK^* \times \KK^*$
act on $X$ by defining the weight of 
the variable $z_i$ as the $i$-th column of the matrix
$$ 
\left[
\begin{array}{rrrrrr}
1 & 1 & 1 & 1 & 1 & 1
\\
1 & 2 & 3 & 3 & 4 & 5
\end{array}
\right]
$$
Note that this action lifts
the action of the second factor 
$\KK^*$ on $G(2,4)$ 
given in $\PP^5$
by homogeneous Pl\"ucker 
Coordinates as
\begin{eqnarray*}
t_2 \mal [z] 
& := & 
[t_2z_1, t_2^2z_2,t_2^3z_3,t_2^3z_4,t_2^4z_5,t_2^5z_6].
\end{eqnarray*} 

The open $T$-invariant 
subsets $U \subseteq X$ 
admitting a complete orbit space
$U/T$ are,
via the tautological projection,
in one-to-one 
correspondence with the 
open $\KK^*$-invariant 
subsets $V \subseteq G(2,4)$
admitting a complete orbit variety
$V/\KK^*$.
The latter ones are well known,
see~\cite{BBSw1}:
there are six of them;
four having a projective 
orbit variety, 
and two having quite exotic
orbit spaces, 
which are not even embeddable
into toric varieties,
see~\cite{Sw1}. 

Let us see how to recover 
this picture via our method.
We need a describing 
pp-divisor fo the $T$-action on 
$X$.
According to the recipe discussed
in~\cite[Section~11]{AlHa}, we first 
determine a pp-divisor for the 
(equivariant) ambient space
$\KK^6$, using the language
of toric varieties. 
We supress the details of 
computation; all of them 
are standard toric geometry,
one may even use, e.g., 
the software
package~\cite{TorDiv} as a 
help.

As the underlying projective 
variety $Y_{\rm ambient}$, we take the normalized
component of the GIT-limit of the 
$T$-action on $\KK^6$.
Concretely $Y_{\rm ambient}$ is the 
toric variety given by the
fan $\Sigma$ 
in $\QQ^4$ having its rays 
through the vectors
$$
v_1 := (1,0,0,0),
\quad
v_2 := (0,1,0,0),
\quad
v_3 := (2,1,0,0),
\quad
v_4 := (3,2,1,1),
$$
$$
a_1 := (-4,-3,-2,-2),
\qquad
a_2 := (0,0,1,0),
\qquad
a_3 := (0,0,0,1). 
$$
The fan $\Sigma$ comprises
twelve maximal cones.
In terms of the above vectors, 
they are given by
$$
\begin{array}{llll}
\cone(v_1,a_2,a_1,v_3), 
& 
\cone(v_2,a_2,a_1,v_3), 
& 
\cone(v_1,a_2,v_4,v_3), 
& 
\cone(v_2,a_2,v_4,v_3), 
\\
\cone(v_1,a_3,a_1,v_3), 
&
\cone(v_2,a_3,a_1,v_3), 
&
\cone(v_1,a_3,v_4,v_3), 
&
\cone(v_2,a_3,v_4,v_3), 
\\
\cone(v_2,a_2,a_3,v_4), 
&
\cone(v_2,a_2,a_3,a_1), 
&
\cone(v_1,a_2,a_3,v_4), 
&
\cone(v_1,a_2,a_3,a_1). 
\end{array}
$$
Denoting by $D_1, \ldots, D_4$
the invariant prime divisors
of $Y_{\rm ambient}$ 
corresponding to the rays 
through $v_1, \ldots, v_4$,
we obtain a describing pp-divisor 
\begin{eqnarray*}
\mathfrak{D}_{\rm ambient}
& = & 
\Delta_1 \otimes D_1 + \ldots \Delta_4 \otimes D_4,
\end{eqnarray*}
where the (common)
tail cone of the
polyhedra $\Delta_i \subset \QQ^2$
is generated by the vectors
$(-1,1)$ and $(5,-1)$ and, thus, they
are given by
$$
\begin{array}{ll}
\vertices(\Delta_1) 
= 
\{(0,0),(2,-1)\},
&
\vertices(\Delta_2) 
= 
\{(-1,1)\},
\\
\vertices(\Delta_3) 
= 
\{(0,0),(3,-1)\},
&
\vertices(\Delta_4) 
= 
\{(0,0),(4,-1)\}.
\end{array}
$$ 

It turns out that
for $\KK^6$ we have 
we have four open subets 
$W_1, \ldots, W_4 \subset \KK^6$ 
with a complete
(in fact projective) orbit 
variety $W/T$.
These arise from the following
four coherent
collections 
(the vertices are 
listed according to the 
enumeration 
$\Delta_1, \ldots, \Delta_4$):
\begin{eqnarray*}
\{(0,0), \; (-1,1), \; (0,0), \; (0,0)\},
& \quad & 
\{(2,-1), \; (-1,1), \; (0,0), \; (0,0)\},
\\
\{(2,-1), \; (-1,1), \; (3,-1), \; (0,0)\},
& \quad &
\{(2,-1), \; (-1,1), \; (3,-1), \; (4,-1)\}.
\end{eqnarray*}

A pp-divisor $\mathfrak{D}$ 
for the $T$-action 
on $X$ lives on the (normal)
closure $Y$ of the 
image of the intersection
$X \cap (\KK^*)^6$ in 
$Y_{\rm ambient}$, and
$\mathfrak{D}$ can be taken 
as the pull back of 
$\mathfrak{D}_{\rm ambient}$
with respect to the 
inclusion 
$\imath \colon Y \to Y_{\rm ambient}$.
Pulling back toric prime
divisors $D_i$ gives
$$ 
\imath^*(D_1) = E_1,
\quad
\imath^*(D_2) = E_2,
\quad
\imath^*(D_3) = E_3^a \cup E_3^b,
\quad
\imath^*(D_4) = E_4,
$$ 
with prime divisors 
$E_1,E_2,E_4$ and 
$E_3^a$, $E_3^b$,
the latter two being
disjoint from each 
other.
The pp-divisor for 
the $T$-action on 
$X$ then is given by
\begin{eqnarray*}
\mathfrak{D}
& = & 
\Delta_1 \otimes E_1 + 
\Delta_2 \otimes E_2 + 
\Delta_3 \otimes E_3^a + 
\Delta_3 \otimes E_3^b + 
\Delta_4 \otimes E_4.
\end{eqnarray*}
Up to the splitting 
of $D_3$ into two 
disjoint components,
the intersection 
behaviour of the pulled
back divisors is as 
before, which can be directly 
checked in toric affine 
charts.
This gives six coherent 
collections of vertices
(listed according to 
the enumeration
$\Delta_1,\Delta_2,\Delta_3^a,\Delta_3^b,\Delta_4$):
$$
\{(0, 0),  (-1, 1),  (0, 0),  (0, 0),  (0, 0)\},
\
\{(2, -1),  (-1, 1),  (0, 0),  (0, 0),  (0, 0)\},
$$
\{(2, -1),  (-1, 1),  (0, 0),  (3, -1),  (0, 0)\},
\
\{(2, -1),  (-1, 1),  (3, -1),  (0, 0),  (0, 0)\},
$$
\{(2, -1),  (-1, 1),  (3, -1),  (3, -1),  (0, 0)\},
\
\{(2, -1),  (-1, 1),  (3, -1),  (3, -1),  (4, -1)\}.
$$
\end{example}

The reader might be a little
disappointed about the computational
efforts needed for the preceeding 
example.
The situation turns much better, if 
one considers actions of tori having 
generic orbits of small 
codimension (instead of small 
dimension):

\begin{example}
Let $X$ be a normal, affine
variety with a good effective 
action of a torus $T$ such that
$\dim(T) = \dim(X) - 1$ holds.
Then $X$ arises from a 
pp-divisor $\mathfrak{D}$
on a projective curve $Y$,
and $\mathfrak{D}$ is of the form
\begin{eqnarray*}
\mathfrak{D}
& = & 
\sum_{i=1}^r \Delta_i \otimes \{y_i\},
\end{eqnarray*}
where the $y_i \in Y$ are 
pairwise different points.
Any collection $v_1, \ldots, v_r$
of  vertices $v_i \in \Delta_i$ 
is coherent, and hence the 
collection 
of open $T$-invariant
$U \subseteq X$ with
complete orbit space $U/T$ 
is in bijection to the set
$$ 
\vertices(\Delta_1)
\times \ldots \times
\vertices(\Delta_r).
$$
\end{example}

\end{document}